\documentclass{article}
\usepackage{amssymb,amsmath,algorithm,algorithmic,enumerate}

\usepackage{natbib} 
\usepackage{algorithm} 
\usepackage{algorithmic} 
\usepackage{multirow}

\title{New Pivot Selection for Sparse Symmetric Indefinite Factorization}
\author{Duangpen Jetpipattanapong, Gun Srijuntongsiri}
\begin{document}

\maketitle

\begin{abstract} 
We propose a new pivot selection technique for symmetric indefinite factorization of sparse matrices.
Such factorization should maintain both sparsity and numerical stability of the factors, both of which depend solely on the choices of the pivots.
Our method is based on the minimum degree algorithm and also considers the stability of the factors at the same time.
Our experiments show that our method produces factors that are sparser than the factors computed by MA57 \cite{MA57} and are stable.
\end{abstract}

\section{Introduction}

There are two general approaches to solving linear systems: direct and iterative methods.
Iterative methods start with an initial guess and successively generate better approximate solutions at each iteration.
The running time of an iterative method depends directly on the required accuracy of the solution.
Direct methods, on the other hand, attempt to find the solution by a finite number of operations and usually involve factoring the matrix.  

Our work here is on direct methods for solving linear systems where the matrix is symmetric indefinite and sparse, which has many applications in linear and nonlinear optimization and finite element computation, for example.  Symmetric indefinite factorization (SIF) is not unique as the resulting factors depend on the choices of the pivots during the factorization.  Pivots should be chosen such that the resulting factors are stable and do not have many fill-ins---the entries that are zeros in the original matrix but are nonzeros in the factors.  

There are many heuristic techniques for selecting pivots to minimize the number of fill-ins for the related problem of Cholesky factorization, which is the most suitable factorization for symmetric positive definite matrices, in literature.  We briefly discuss a few such well-known techniques here since some of their ideas are also applicable to SIF.
These ordering algorithms can be classified into three classes: local, global, and hybrid approaches.
Local approach such as the minimum degree and the minimum fill algorithms \cite{MD1, MD2, MD3, MDevol, MFNP} selects the pivot that is expected to minimize the number of fill-ins at each factorization step in a greedy fashion.  Global approach such as Cuthill-McKee and nested dissection methods \cite{CM, ND, NDgen} selects pivots by considering the overall structure of the matrix.  Hybrid approach, on the other hand, combines the ideas from both local and global approaches.

The well-known minimum degree algorithm \cite{MD1} chooses the column that has the minimum off-diagonal nonzero elements in the remaining matrix as the pivot for the current step.  
Different improvements of the minimum degree algorithm have been proposed \cite{MDevol} such as multiple minimum degree \cite{MMD1} and approximate minimum degree algorithms \cite{AMD1} and become the practical standard in the implementations.

Another famous pivot selection algorithm is nested dissection \cite{ND}.  By defining a graph whose vertices represent each column of the matrix and whose edges represent nonzero entries in the matrix, nested dissection recursively find a separator---a set of vertices that partitions the graph into two disconnected subgraphs---and ordering the pivots recursively with the two subgraphs first followed by the separator vertices.  Cuthill-McKee \cite{CM} propose another pivot selection algorithm that aims to reduce the bandwidth of the matrix based on breadth first search of the  structure graph.

The main difference between Cholesky factorization and SIF is in the size of pivots.  For SIF, each pivot can be either a scalar or a 2-by-2 matrix while pivots in Cholesky factorization are all scalar.
Moreover, unlike Cholesky factorization, the choice of pivots in SIF also affects the stability of the resulting factors \cite{SIF1}.

There are many pivot selection algorithms proposed specifically for SIF such as 
Bunch-Parlett \cite{BP}, Bunch-Kaufman \cite{BK}, and bounded Bunch Kaufman (BBK) \cite{BBK} algorithms.
Bunch-Parlett method searches the whole remaining submatrix at each stage for the largest-magnitude diagonal and the largest-magnitude off-diagonal.  It chooses the largest-magnitude diagonal as the 1-by-1 pivot if the resulting growth rate is acceptable.  Otherwise, it selects the largest-magnitude off-diagonal and its relative diagonal elements as the 2-by-2 pivot block. 
This method requires $O(n^3)$ comparisons and yields a matrix  $L$ whose maximum element is bounded by 2.781.
Bunch-Kaufman pivoting strategy searches for the largest-magnitude off-diagonal elements of at most two columns for each iteration. It requires $O(n^2)$ comparisons, but the elements in $L$ are unbounded. 
BBK combines the two above strategies.  By monitoring the size of the elements in $L$, BBK uses the Bunch-Kaufman strategy when it yields modest element growth. Otherwise, it repeatedly searches for an acceptable pivot.
In average cases, the total cost of BBK is the same as Bunch-Kaufman, but in the worst cases its cost can be the same as that of the Bunch-Parlett strategy.  

Additionally, there are other types of techniques for solving symmetric indefinite linear systems.  Paige and Saunders \cite{SSIF1} propose two algorithms, SYMMLQ and MINRES, for solving such systems.  The algorithms apply orthogonal factorization together with the conjugate gradient method to solve the system.  Duff et al.\ \cite{SSIF2} propose a pivotal strategy for decomposing sparse symmetric indefinite matrices that limits the magnitude of the element in the factors for stability.
Duff and Reid \cite{MFrontal1} propose a multifrontal method to solve indefinite sparse symmetric linear systems based on minimum-degree ordering.
The multifrontal approach is widely used in many sparse direct solvers, for example, MA57 and MUMPS \cite{DS1, MA57}.

We propose a new pivot selection algorithm for sparse SIF.  Our algorithm applies the idea of minimum-degree ordering to consider both 1-by-1 and 2-by-2 pivots while also considers the stability of the resulting factors.  Our experiments show that our algorithm produces stable factors that are sparser than the factors produced by \cite{MA57}.

For the rest of the article, we describe symmetric indefinite factorization in Section \ref{sec2}.  Section \ref{sec3} explains the minimum-degree ordering algorithm.  Section \ref{sec4} describes our algorithm.  Section \ref{sec5} shows our experiment and the results.  Finally, we conclude the article in Section \ref{sec6}.

\section{Symmetric Indefinite Factorization}\label{sec2}
Solving a linear system
\begin{equation}\label{eq1}
Ax = b,
\end{equation}
where $A \in \mathbb{R}^{n \times n}$ is symmetric indefinite, is generally done by first obtaining the symmetric indefinite factorization 
\begin{equation}\label{eq2}
P^{T}AP = LBL^{T},
\end{equation}
where $P$ is a permutation matrix, $L$ is a unit lower triangular matrix, and $B$ is a block diagonal matrix
\[
B=
\begin{bmatrix} B^{(1)}      & 0          &\cdots  & 0        \\     
                     0          & B^{(2)}      &\cdots  & 0        \\
                     \vdots  & \vdots &\ddots  & \vdots\\
                     0          & 0          &\cdots  & B^{(K)}    \end{bmatrix},
\]
where each block $B^{(k)}$ is either a 1-by-1 or 2-by-2 matrix and is nonsingular.  Matrix $P$ represents the pivoting---the exchanging of rows and columns of $A$---during the factorization in order to maintain both sparsity and stability of the factor.
After obtaining the factorization, back and forward substitutions are used to compute the solution of (\ref{eq1}) by the following steps:
\begin{enumerate}
\item[(i)]\label{step1} Solve $z$ : $Lz  =   P^T b$.
\item[(ii)]\label{step2} Solve $\hat{z}$ : $B\hat{z} = z$.
\item[(iii)]\label{step3}Solve $\bar{z}$ : $L^T\bar{z}= \hat{z}$.
\item[(iv)]\label{step4}Set : $x=P\bar{z}$.
\end{enumerate}

Recall that Steps (ii) and (iv) are trivial (due to the structure of $P$ and $B$) and therefore the computational time for solving the linear system depends solely on the factorization and back and forward substitutions in Steps (i) and (iii), which in turn depend on the sparsity of $L$.  

To perform symmetric indefinite factorization, let $A^{(k)}$ be the (smaller) matrix that remains to be factorized in the $k$th iteration. The algorithm starts with $A^{(1)}=A$. 
For each iteration, we first identify a submatrix $B^{(k)}$ from elements of $A^{(k)}$ that is suitable to be used as the pivot block.
The submatrix $B^{(k)}$ is either a single diagonal element of $A^{(k)}$  
$\left(\left\lbrack a^{(k)}_{ll} \right\rbrack\right)$
 or a 2-by-2 block with two diagonal elements of $A^{(k)}$  
$\left(\begin{bmatrix} a^{(k)}_{ll}  & a^{(k)}_{lr} \\  a^{(k)}_{rl} & a^{(k)}_{rr}  \end{bmatrix}\right)$. 
Note that there are many methods to select the pivot (We explain pivot selection algorithms and our proposed pivot selection algorithm in the next section).  Next, we find the permutation matrix $P^{(k)}$ satisfying
\begin{equation} \label{eq3}
(P^{(k)})^TA^{(k)}P^{(k)}=\begin{bmatrix} B^{(k)}     & (C^{(k)})^T   \\      C^{(k)}  & Z^{(k)}       \end{bmatrix}.
\end{equation}
That is, $P^{(k)}$ is the permutation matrix corresponding to the exchanges of rows and columns that move the chosen pivot to the top-left corner.  The right-hand side of (\ref{eq3}) can be factorized as
\begin{eqnarray} \label{eq4}
(P^{(k)})^TA^{(k)}P^{(k)}&=&\begin{bmatrix} I  & 0   \\  C^{(k)}(B^{(k)})^{-1}    & I \end{bmatrix}\cdot
\begin{bmatrix} B^{(k)}  & 0   \\      0    & Z^{(k)}-C^{(k)}(B^{(k)})^{-1} (C^{(k)})^T  \end{bmatrix}\cdot\nonumber\\
&&\begin{bmatrix} I  & (B^{(k)})^{-1}(C^{(k)})^T   \\   0     & I  \end{bmatrix}.
\end{eqnarray}
Let $L^{(k)} = C^{(k)}(B^{(k)})^{-1}$ and $A^{(k+1)}$ = $Z^{(k)}-C^{(k)}(B^{(k)})^{-1}(C^{(k)})^T$. Equation (\ref{eq4}) can be rewritten as
\begin{equation} \label{eq5}
(P^{(k)})^TA^{(k)}P^{(k)}=
\begin{bmatrix} I  & 0   \\     L^{(k)}   & I \end{bmatrix}\cdot
\begin{bmatrix} B^{(k)}  & 0   \\      0    & A^{(k+1)} \end{bmatrix}\cdot
\begin{bmatrix} I  &    (L^{(k)})^T \\0& I \end{bmatrix}.
\end{equation}
The same process can be repeated recursively on the matrix $A^{(k+1)}$. Note that the dimension of  $A^{(k+1)}$ is less than the dimension of $A^{(k)}$ by either one or two depending on the dimension of $B^{(k)}$.

\section{Pivot selection with Minimum Degree} \label{sec3}
Finding the optimal ordering that minimizes fill-in is NP-hard \cite{MFNP} therefore a heuristic is often used for pivot selection. 
Choosing pivot at each step should be inexpensive, lead to at most modest growth in the elements of the remaining matrix, and not cause $L$ to be too much denser than the original matrix.  One of the well-known and efficient pivot selection techniques is the minimum degree algorithm \cite{MD1, MD2, MD3}.
The algorithm considers the pivot based on the following graph model. 
Define an undirected graph $G=(V,E)$, where $V=\{1,...,n\}$ and $E=\{\{i,j\}:i\neq j$ and $a_{ij}\neq 0\}$.
Observe that the degree of $v$ ($\text{deg}(v)$), where $v \in V$, is the number of nonzero off-diagonal elements on the $v$th row.
The vertex $v$ with minimum $\text{deg}(v)$ is chosen as the pivot.

Define the elimination graph $G_v=(V\setminus \{v\},E')$, where $E'=E\cup\{\{i,j\}:\{i,v\}\in E \mbox{ and }\{v,j\}\in E\}\setminus \{\{v,i\}:i=1,2,...,n\}$.
Graph $G_v$ is used to choose the next pivot, and so on.
That is, the minimum degree algorithm is as follows. 

\begin{algorithm}
\caption{Minimum Degree Algorithm}\label{algo1}
\begin{algorithmic}
\STATE Define $G$ as described above.
\WHILE{ $G \neq \emptyset$}
\STATE  $v$ = the vertex with minimum $\text{deg}(v)$
 \STATE  $G=G_v$
\ENDWHILE
\end{algorithmic}
\end{algorithm}
Note that the minimum degree algorithm identifies the pivot at each step without any numerical calculation.
For this reason, it can be used as the ordering step before factorizing the matrix.

Many improvements of the minimum degree algorithm and its implementation have been proposed \cite{MDevol} such as decreasing the computation time for the degree update by considering the indistinguishable nodes \cite{massEli} or minimum degree independent nodes \cite{MMD1}, reducing the computation cost by using an approximate minimum degree \cite{AMD1}, and saving space by using the quotient graph model \cite{quotGraph}.

\section{Our pivot selection algorithm} \label{sec4}

Unlike in Cholesky factorization, pivots in symmetric indefinite factorization can be either a scalar or a $2$-by-$2$ matrix therefore the minimum degree algorithm cannot be used as is in this case. 


The stability condition that our algorithm uses is proposed by Duff et~al.\ \cite{stCond} and also used as a thresholding test for 1-by-1 and 2-by-2 pivots in MA57 \cite{MA57}.
We consider a 1-by-1 pivot $a_{ii}$ to be \emph{acceptably stable} if
\begin{equation} \label{eq6}
|a_{ii}|\geq \alpha\max_{r\neq i}  |a_{ri}|.
\end{equation}
Similarly, a 2-by-2 pivot $\begin{bmatrix} a_{ii}&a_{ij}\\a_{ji}&a_{jj}\end{bmatrix}$ is considered to be \emph{acceptably stable} if 
\begin{equation}\label{eq7}
\left|\begin{bmatrix} a_{ii}&a_{ij}\\a_{ji}&a_{jj}\end{bmatrix}^{-1}\right|
\cdot \begin{bmatrix}\max_{r\neq i, r\neq j} |a_{ri}|\\ \max_{r\neq i, r\neq j}|a_{rj}|\end{bmatrix}
\leq
\begin{bmatrix}\alpha^{-1}\\\alpha^{-1}\end{bmatrix}.
\end{equation}
Conditions (\ref{eq6}) and (\ref{eq7}) limit the magnitudes of the entries of $L$ to $1/\alpha$ at most.
The appropriate value of $\alpha$ is $0<\alpha\leq 0.5$.
The default value of $\alpha$ in MA57 is 0.01 \cite{MA57}.

Let us call the column with the fewest number of off-diagonal nonzeros the \textit{minimum degree column}.
Let $i$ be the minimum degree column of the matrix $A$.
We accept $a_{ii}$ as the 1-by-1 pivot $(B^{(k)})$ if $a_{ii}$ satisfies (\ref{eq6}).
Otherwise, we proceed to search for a suitable 2-by-2 pivot $\begin{bmatrix} a_{ii}  & a_{ij} \\  a_{ji} & a_{jj}\end{bmatrix}$ that satisfies (\ref{eq7}) as follows.
Let 
\begin{equation}
\label{eq10}
Z_i = \{ z | a_{iz}\neq 0 \mbox{ and } z\neq i\}.
\end{equation}
Consider all submatrices $\begin{bmatrix} a_{ii}&a_{iz}\\a_{zi}&a_{zz}\end{bmatrix}$, where $z \in Z_i$, as the candidates for a $2$-by-$2$ pivot.  The degree of each candidate $\text{deg}(i,z)$ is the number of rows $l$ where $l \neq i, z$ and at lease one of $a_{li}$ and $a_{lz}$ is nonzero.  To compute $\text{deg}(i,z)$, define
\begin{equation} \label{eq8}
d(i,z,l)=\left\{
\begin{array}{ll}
0 ,& \mbox{if } a_{li}=0 \mbox{ and } a_{lz} = 0,\\
1 ,& \mbox{otherwise.} 
\end{array}
\right.
\end{equation}
Hence,
\begin{equation} \label{eq9}
\text{deg}(i,z) = \sum_{l\neq i,z}d(i,z,l).
\end{equation}
Our algorithm then considers all of the candidates with the minimum out-degree.  Specifically, 
$\begin{bmatrix} a_{ii}  & a_{ij} \\  a_{ji} & a_{jj}\end{bmatrix}$ is qualified if
\begin{equation}\label{eq11}
\text{deg}(i,j) = \min_{z \in Z_i} \text{deg}(i,z).
\end{equation}
If a qualified candidate also satisfies (\ref{eq7}), it is chosen as a pivot.  Otherwise, we remove $j$ from the $Z_i$ and repeat the process of selecting a $2$-by-$2$ pivot until we either find a qualified candidate that also satisfies (\ref{eq7}) or $Z_i$ becomes empty.  In the latter case, we set $i$ to be the next next minimum degree column and repeat the process from the beginning (from testing whether $a_{ii}$ is a suitable $1$-by-$1$ pivot).  
The algorithm is as shown in Algorithm \ref{algo2} below.

\begin{algorithm}
\caption{Our Pivot Selection Algorithm}\label{algo2}
\begin{algorithmic}
\STATE // $A$ is a $n$-by-$n$ symmetric indefinite matrix
\STATE Let $M = \{1, 2, ... , n\}$
\WHILE{a suitable pivot is not yet found and $M$ is not empty}
\STATE Let $i$ be the minimum degree column among all column indices in $M$
\IF{$a_{ii}$ is accepted}
   \STATE Use $a_{ii}$ as the 1-by-1 pivot
\ELSE
   \STATE Let $Z_i = \left\{ z | a_{iz}\neq 0 \mbox{ and } z\neq i\right\}$ 
   \WHILE {a suitable pivot is not yet found and $Z_i$ is not empty}
      \STATE Let $j$ be such that $\begin{bmatrix} a_{ii}  & a_{ij} \\  a_{ji} & a_{jj}\end{bmatrix}$ has the minimum out-degree and $j \in Z_i$
      \IF{$\begin{bmatrix} a_{ii}  & a_{ij} \\  a_{ji} & a_{jj}\end{bmatrix}$ satisfies (\ref{eq7})}
         \STATE Use $\begin{bmatrix} a_{ii}  & a_{ij} \\  a_{ji} & a_{jj}\end{bmatrix}$ as the 2-by-2 pivot
      \ELSE
         \STATE Remove $j$ from $Z_i$ 
      \ENDIF
   \ENDWHILE
   \STATE Remove $i$ from $M$
\ENDIF
\ENDWHILE
\end{algorithmic}
\end{algorithm}
Lastly, when the remaining matrix is fully dense, we continue with a conventional pivot selection algorithm such as BBK instead.


\section{Experiments and results} \label{sec5}

This section compares the efficiency of our algorithm with MA57, which is based on the multifrontal method. 
The experiments are performed in Matlab 2011a on matrices of varying dimensions from 100 to 5000.
For each dimension, we vary the percentage of nonzeros in the matrices from 5 to 30 percent.
We test with 20 different instances for problems with 100, 300, and 500 dimensions and 10 different instances for problems with 1000,3000, and 5000 dimensions.  We show the percentage of nonzeros in the factor $L$ of the two methods in Table \ref{table1}, which shows that our method produces sparser factors than MA57 in all cases.  Note that the small percentage improvement for large matrices are not insignificant as small decrease in nonzeros does lead to significantly faster factorization time.  Finally, Table \ref{table3} shows the residuals $\left\| P^TAP - LBL^T\right\|$ (* Is this how you compute residuals?  If not, change to the one you use.*) of the results of both methods. The result shows that our method produces more accurate factors than MA57.

\begin{table}
\caption{Average percentage of nonzeros in the factor $L$ produced by MA57 and our algorithm for problems with 100, 300, 500, 1000, 3000, and 5000 dimensions and 30, 20, 10, and 5 percent of nonzeros in the matrix.  The percentage of nonzeros in $L$ is computed by dividing the number of nonzeros in $L$ by $n^2$ and then multiplying the result by $100$.}\label{table1}
\begin{tabular}{|c|r|r|r|r|r|r|r|r|}
\hline
\multirow{3}{*}{$n$} 
&\multicolumn{8}{c|}{Percentage of nonzeros in $L$}\\
\cline{2-9}
&\multicolumn{2}{c|}{30}
&\multicolumn{2}{c|}{20}
&\multicolumn{2}{c|}{10}
&\multicolumn{2}{c|}{5}\\
\cline{2-9}
&MA57
&Our method
&MA57
&Our method
&MA57
&Our method
&MA57
&Our method\\
\hline
100&46.20&45.54&40.90&39.24&22.68&18.73&11.02&6.60\\
300&46.07&45.39&43.03&41.89&35.76&33.15&25.17&21.23\\
500&47.37&46.98&45.26&44.52&39.98&38.17&17.42&12.04\\
1000&48.53&48.35&47.39&47.00&44.01&43.02&38.46&36.36\\
3000&49.46&49.37&49.00&48.84&47.56&47.21&45.11&44.19\\
5000&49.64&49.61&49.36&49.26&48.47&48.22&46.86&46.23\\
\hline
\end{tabular}
\end{table}

\ifx 
Table 2 show the factorization time of both two methods.
The result show that, when the problem is more sparse our algorithm is save more time comparing with MA57. 
\begin{table}
\caption{Average factorization time of MA57 and our algorithm for problems with 300, 500, 1000, and 2000 variables and 70, 80, 90 and 95 percent of zeros in the matrix.}
\begin{tabular}{|c|r|r|r|r|r|r|r|r|}
\hline
\multirow{3}{*}{$n$} 
&\multicolumn{8}{c|}{factorization time (s)}\\
\cline{2-9}
&\multicolumn{4}{c|}{MA57}
&\multicolumn{4}{c|}{Our method}\\
\cline{2-9}
&70
&80
&90
&95
&70
&80
&90
&95\\
\hline
100	&	0.00059	&	0.00056	&	0.00038	&	0.00029	&	0.01631	&	0.02077	&	0.02182	&	0.02176	\\													
300	&	0.00523	&	0.00499	&	0.00432	&	0.00275	&	0.11447	&	0.13355	&	0.14325	&	0.13120	\\													
500	&	0.01507	&	0.01547	&	0.01283	&	0.00579	&	0.48870	&	0.57355	&	0.66205	&	0.40296	\\													
1000	&	0.09170	&	0.08992	&	0.07798	&	0.06432	&	3.26590	&	3.99701	&	5.02232	&	5.48706	\\													
3000	&	1.84638	&	1.92231	&	1.87029	&	1.76557	&	47.71865	&	62.89543	&	96.95069	&	120.06834	\\													
5000	&	6.56945	&	6.68904	&	6.62526	&	6.41201	&	159.18106	&	201.59030	&	305.14178	&	409.57366	\\													\hline
\end{tabular}
\end{table}
\fi

\begin{table}
\caption{Average residuals of the factorization produced by MA57 and our algorithm for problems with 300, 500, 1000, and 2000 dimensions and 30, 20, 10 and 5 percent of nonzeros in the matrix.} \label{table3}
\begin{tabular}{|c|r|r|r|r|r|r|r|r|}
\hline
\multirow{3}{*}{$n$} 
&\multicolumn{8}{c|}{Residual ($\times 10^{-10}$)}\\
\cline{2-9}
&\multicolumn{2}{c|}{30}
&\multicolumn{2}{c|}{20}
&\multicolumn{2}{c|}{10}
&\multicolumn{2}{c|}{5}\\
\cline{2-9}
&MA57
&Our method
&MA57
&Our method
&MA57
&Our method
&MA57
&Our method\\
\hline
100	&	0.00339&	0.00018	&	0.00410&	0.00022	&	0.00190	&	0.00016&	0.00045	&	0.00006	\\													
300	&	0.03072&	0.00077	&	0.02634&	0.00083	&	0.02378&	0.00083	&	0.01039	&	0.00059	\\													
500	&	0.08489&	0.00128	&	0.06665&	0.00161	&	0.04665&	0.00169	&	0.02199&	0.00076	\\													
1000	&	0.20679&	0.00342	&	0.21691&	0.00374	&	0.17399&	0.00333	&	0.10509&	0.00355	\\													
3000	&	1.63656	&	0.01312&	1.80491	&	0.01281&	1.32961	&	0.02150&	1.13003	&	0.02160	\\													
5000	&	4.45974&	0.02488	&	3.49949&	0.02361	&	2.51524&	0.03264	&	2.20916&	0.03152	\\													\hline
\end{tabular}
\end{table}

\ifx
\begin{table}
\caption{Average factorization time of our method using Algorithm 2 and Algorithm 3 for problems with 300, 500, 1000, and 2000 variables, 95 percent of zeros in the matrix and 70, 80, 90 and 95 percent of zeros in the matrix.}
\begin{tabular}{|c|r|r|r|r|r|r|r|r|}
\hline
\multirow{3}{*}{$n$} 
&\multicolumn{8}{c|}{factorization time (s)}\\
\cline{2-9}
&\multicolumn{4}{c|}{Algorithm1}
&\multicolumn{4}{c|}{Algorithm2}\\
\cline{2-9}
&70
&80
&90
&95
&70
&80
&90
&95\\
\hline
100	&	0.01631	&	0.02077	&	0.02182	&	0.02176	&	0.02800	&	0.03388	&	0.04940	&	0.05190	\\													
300	&	0.11447	&	0.13355	&	0.14325	&	0.13120	&	0.19929	&	0.21520	&	0.26964	&	0.26227	\\													
500	&	0.48870	&	0.57355	&	0.66205	&	0.40296	&	0.81153	&	0.78094	&	0.83080	&	0.73690	\\													
1000	&	3.26590	&	3.99701	&	5.02232	&	5.48706	&	3.56344	&	4.33894	&	5.78310	&	6.46104	\\													
3000	&	47.71865	&	62.89543	&	96.95069	&	120.06834	&	234.17849	&	66.50211	&	99.11141	&	137.02838	\\													
5000	&	159.18106	&	201.59030	&	305.14178	&	409.57366	&	177.53116	&	217.69066	&	388.44482	&	410.39756	\\													\hline
\end{tabular}
\end{table}
\fi
\ifx
For the problems which require more accuracy. The pivot selection needs to compare in many times. 
Algorithm 3 is our improve algorithm with suitable to the problem which require many accuracy. 
Table 3 shows the factorization time of the our method and our improve method in Algorithm 3.
We compare the problem with 95 percent sparse by varying the size of variable with $\alpha$ is equal to 0.01, 0.1, 1, and 10 respectively. 
The result shows that our improved algorithm is saving more time when the conditional value of $\alpha$ is larger. 
\fi

\section{Conclusion} \label{sec6}

In this article, we propose a new pivot selection algorithm for symmetric indefinite factorization. Our method is based on the minimum degree algorithm but is able to select both 1-by-1 and 2-by-2 pivots that are stable. Our experimental results show that our algorithm produces factors that are stable and also sparser than MA57.

\bibliography{myrefs}
\end{document}